\documentclass[11pt]{amsart}

\usepackage{amsmath, amsthm, amssymb, amsfonts, enumerate,stmaryrd}
\usepackage[colorlinks=true,linkcolor=blue,urlcolor=blue]{hyperref}
\usepackage{dsfont}
\usepackage{color}
\usepackage{geometry}
\usepackage{todonotes}
\usepackage{soul}
\usepackage{graphicx,epstopdf}
\usepackage{caption}
\usepackage{mathtools}
\usepackage{bbm}

\usepackage{appendix}

\usepackage{tikz}
\tikzset{every node/.style={align=center}}

\mathtoolsset{showonlyrefs}

\makeatletter
\@namedef{subjclassname@2020}{%
  \textup{2020} Mathematics Subject Classification}
\makeatother

\newcommand\independent{\protect\mathpalette{\protect\independenT}{\perp}}
\def\independenT#1#2{\mathrel{\rlap{$#1#2$}\mkern2mu{#1#2}}}

\def \ud{\mathrm{d}}

\newtheorem{theorem}{Theorem}[section]
\newtheorem{lemma}[theorem]{Lemma}
\newtheorem*{lemma*}{Lemma}

\newtheorem{proposition}[theorem]{Proposition}
\newtheorem{definition}[theorem]{Definition}
\newtheorem{remark}[theorem]{Remark}
\newtheorem{assumptions}[theorem]{Assumption}

\title[Reflected Elastic diffusions]{reflected Elastic diffusions}
\author[M.\ Allemandri]{Mattia Allemandri}
\address{M.\ Allemandri: School of Management and Economics, Dept.~ESOMAS, University of Turin, C.so Unione Sovietica 218bis, 10134 Turin, Italy.}
\email{\href{mailto:mattia.allemandri@unito.it}{mattia.allemandri@unito.it}}

\numberwithin{equation}{section}

\begin{document}

\begin{abstract}
We study the reflected elastic diffusion (RED), the pathwise counterpart of Feller's elastic boundary condition for a one-dimensional diffusion on $[0,\infty)$ with time-homogeneous, continuous drift and diffusion coefficients.
Given a 
diffusion $X$ reflecting at zero, we construct a Feller process $\widetilde X$ by killing $X$ at the first time its 
local time at zero exceeds an independent exponential random variable. We characterize the  infinitesimal generator of $\widetilde X$ explicitly and prove that its domain satisfies a Robin boundary condition,  thereby identifying $\widetilde X$ as a probabilistic representation of Feller's analytical characterization of the elastic boundary. 
A key analytical ingredient, of independent interest, is a pair of moment asymptotics, as $t\downarrow 0$, for the boundary local time of a 
reflected diffusion, extending known results for the unreflected case.
Uniqueness in law of $\widetilde X$ among Feller REDs sharing the same interior dynamics is established through a resolvent/boundary-value-problem argument. The result extends the classical Itô--McKean construction of elastic Brownian motion to time-homogeneous diffusions, closing the gap between analytic and probabilistic treatments of elastic killing used in recent applied models.
\end{abstract}

\maketitle

\section{Introduction}\label{sect: intro}
The idea of prescribing the behaviour of a Markov process at the boundary of
its state space through an analytic condition on the domain of its
infinitesimal generator, and then recovering a pathwise realization of that
behaviour, goes back to the foundational work of Feller on one-dimensional diffusions in the
early 1950s \cite{feller1952parabolic, feller1954diffusion,
feller1957generalized}, later systematically developed from a probabilistic perspective by It\^o and McKean \cite{ito1963brownian, ito2012diffusion}. Given an increasing continuous function $S$ (the
scale function) and an increasing function $m$ (the speed measure) on an
interval $(l,r)$, $-\infty \le l < r \le \infty$, Feller obtained a full
classification of the nature of the 
boundary points\footnote{Later on, when
moving from a purely analytical to a probabilistic understanding of this
classification, these boundary points will be referred to as 
the ones of the one-dimensional diffusion characterized via $S$ and $m$.}
$l$ and $r$ 
by asking the question of when the differential operator $\mathcal{G}:= \frac{d}{dm}\frac{d}{dS}$ is the
\emph{infinitesimal generator} of a (strongly continuous,
positivity-preserving, contraction) semigroup $(P_t)_{t\ge0}$ (see
\cite{Mandl1968AnalyticalTO}), i.e., when is it that $\mathcal{G}$ can be
obtained as
$$ \lim_{t\downarrow 0} t^{-1}(P_t f - f) =: \mathcal{L} f $$
for some semigroup $(P_t)_{t\ge0}$ of the type just introduced. When
$l$ is a regular boundary point in the classification induced by $S$ and
$m$ (analogous results hold for $r$), it was shown --- see \cite{Mandl1968AnalyticalTO} for a systematic
account --- that $\mathcal{G}f = \mathcal{L}f$ holds (for some $(P_t)_{t\ge 0}$) on the class
$D(\mathcal{L})$ consisting of those functions $f$ for which $\mathcal{L}f$
exists in some continuous class, and that in addition satisfy
the \emph{(local) boundary condition}\footnote{ A more general boundary condition of non-local type, i.e., involving evaluations of $f$ away from $l$, can also be considered (cf.~\cite[Sec.~3.4]{Mandl1968AnalyticalTO}). Such behaviors, however, lie outside the scope of this paper.}(cf.~\cite[Secs.~2.5 and
3.1]{Mandl1968AnalyticalTO})
$$
\kappa f(l) - \delta \frac{df}{dS}^+(l+) + \beta \mathcal{L}f(l)
=  0, \qquad \kappa, \delta, \beta\ge 0 ,\quad  \kappa+\delta + \beta> 0$$
where $$  \quad  \frac{df} {dS}^+(x) :=  \lim_{h\downarrow 0} \frac{f(x+h)-f(x)}{S(x+h) - S(x)}, ~ x\in(l,r) \quad \text{ and } \quad \frac{df}{dS}^+ (l+) := \lim_{x\downarrow l} \frac{df} {dS}^+(x).$$
The terminology attached to $l$ depends on the specific values taken
by the tuple $(\kappa,\delta,\beta)$. Specifically, an \emph{elastic}
boundary corresponds to $\kappa,\delta>0$, $\beta=0$, an
intermediate behaviour between absorption ($\kappa>0, \delta= \beta=0$) and reflection
($\delta>0, \kappa=\beta = 0$).

The companion probabilistic question --- how to construct, path by path, a
process whose generator satisfies a prescribed boundary condition --- was
answered by It\^o and McKean in their complete treatment of Brownian motion
on the half-line for arbitrary $(\kappa,\delta,\beta)$
\cite{ito1963brownian, ito2012diffusion}. Their
construction includes, in particular, the sticky Brownian motion --- a case
for which, according to the historical account in \cite{peskir2015boundary},
Feller himself seems to have had no conjecture --- and the elastic Brownian
motion, obtained by It\^o and McKean by formalizing a private suggestion of
Feller's (again, see \cite{peskir2015boundary}) as a reflected Brownian
motion killed at a rate proportional to its local time at the origin.

In recent years there has been a renewed interest in these two boundary
behaviours. The literature on sticky boundaries has grown through a series
of papers addressing the problem of representing sticky processes as
solutions of SDEs: this question was revisited for the sticky Brownian
motion in \cite{engelbert2014stochastic}, and solved for the first time for
sticky Bessel processes and sticky Feller diffusions in
\cite{peskir2022sticky} and \cite{Peskir2023StickyFD} respectively; both
papers also provide a closed-form expression for the transition density
function. We also mention \cite{salins2017markov}, where a pathwise
representation is studied for a general diffusion process showing both
stickiness and skewness at a (possibly interior) point.

The elastic boundary condition has likewise attracted renewed attention: \cite{dovidio2023elastic} has extended the elastic Brownian motion to the case with (constant) drift, while \cite{colantoni2025elastic} considered an elastic Brownian motion on a $C^2$ domain which, instead of being killed at the boundary, restarts from
a random position in the interior.

Beyond the classical Brownian setting, elastic killing has recently
resurfaced in more general 
systems 
 with non-constant
coefficients, arising from applications rather than from an 
interest in the boundary behaviour itself. In particular, elastic
processes have been used as a modelling device in the following papers
(potentially among others): \cite{hambly2025spde}, whose focus is on
the empirical measure of a system of interacting particles representing
defaultable entities in a credit risk scenario, construct a system of
correlated diffusions with space- and time-dependent drift and diffusion
coefficients, each subject to elastic killing at a fixed rate $\kappa$; instead, 
\cite{fausti2022interacting} consider a related but more elaborate
mechanism, motivated by a model for an epidemic advancing through a
population, in which the killing rate is itself stochastic and the
process reflects along a moving boundary. 
However, the applied papers above refer to such processes as
\emph{elastic} purely on the basis of an analogy with the elastic
Brownian motion construction, without exploring the domain of the
associated generator, and hence without relating their probabilistic
construction to Feller's analytic classification of boundary points.
To the best of our knowledge, the only result in this direction is
provided by \cite{chen2015one}: starting from a Markov process
killed upon reaching a (boundary) point $a$, in Section~4 they construct an
extension of it by piecing together the excursions from $a$ that return
to $a$, generated by a Poisson point process, and then joining the first excursion that either never returns to $a$ or ends in the cemetery state; the
characteristic measure of the Poisson process carries an atom at the
cemetery path, which accounts for the killing rate. They then study the
$L^2$-infinitesimal generator of the resulting process, whose domain is
determined by a flux condition at $a$ which, for Brownian motion on the
half-line, is shown to coincide with the classical Robin condition. This
construction is abstract and excursion-theoretic, and the boundary
condition is explicitly identified as the Robin one only in the Brownian case.

Rather than relying on excursion theory, we propose a constructive approach to directly build elastically killed reflected diffusions with time-homogeneous, continuous coefficients, where the label \emph{elastic} is rigorously justified by identifying the Robin boundary condition in the domain of the generator

This is done as follows:  in Section \ref{sect: setting_main_result} we give precise definition   of the $(\mu, \sigma)$-reflected elastic diffusion and construct a candidate for it by killing a reflected diffusion via an exponential (multiplicative) functional of its boundary local time. In Section \ref{sect: tildeX_properties} we prove it is strong Markov, introduce its transition semigroup and show it is Feller. In Section \ref{sect: local_time_asymptotics} we provide some results on the asymptotic behaviour of the first and the second moment of the boundary local time for the reflected diffusion, which are then exploited in Section \ref{sect: proof_main_thm} to prove that our pathwise construction provides the unique Feller diffusion on $[0,\infty)$ realizing the elastic
condition $\frac{df}{dS}^+ (0+) = k f(0)$ with the prescribed interior dynamics.

\section{Setting and main result}\label{sect: setting_main_result}

In this section, we introduce the notion of \emph{reflected elastic diffusion} (RED) 
and construct a candidate pathwise representation for it. We begin with some preliminaries on 
Feller processes.

Given a locally compact space $E$ with a countable basis (LCCB), we denote by $C_0(E)$ the space of continuous functions 
vanishing at infinity, equipped with the sup-norm $\|\cdot\|_\infty$. A \emph{Feller semigroup} $(P_t)_{t\ge 0}$ on $C_0(E)$ is a strongly continuous contraction semigroup of positive linear operators mapping 
$C_0(E)$ 
into itself (see, e.g.,~\cite[Ch.~3.2]{revuz2013continuous}), and each Feller semigroup is known (cf.~\cite[Ch.~3, Prop.~2.2]{revuz2013continuous}) to uniquely 
determine a homogeneous Markov transition function, called the \emph{Feller transition 
function}. A Markov process 
admitting such a transition function is called a \emph{Feller process}.

The \emph{infinitesimal generator} of a Feller process is given by
the pair (see, e.g., \cite[Ch.~7, Def.~1.1]{revuz2013continuous})
 \begin{align}
    D(\mathcal{L}) &= \Big\{ f \in C_0(E): \exists g\in C_0(E) \text{ s.t. } \lim_{t\downarrow 0}\Big|\Big| \frac{P_t f - f}{t} - g\Big|\Big|_\infty = 0 \Big\} \label{eq: dom_inf_gen_target}, \\
    \mathcal{L} f &= \lim_{t\downarrow 0}\frac{P_t f - f}{t}, \quad \forall f \in D(\mathcal{L}) \label{eq: inf_gen_target}.
\end{align}

Depending on the behaviour the process shows at the boundaries of its state space 
(briefly, \emph{boundary behaviour}), functions in $D(\mathcal{L})$ satisfy precise 
local conditions at the boundary (we refer to \cite[Ch.~2]{Borodin} for a modern discussion 
on the topic). 
We highlight that the approach adopted in definitions \eqref{eq: inf_gen_target}-\eqref{eq: dom_inf_gen_target} takes a given semigroup as its starting point to define the generator. Hence, unlike Feller's original approach for the classification of boundary behaviours (see Section \ref{sect: intro}), no boundary condition needs to be imposed on $D(\mathcal{L})$ for $\mathcal{L}$ to be well defined as a generator. Rather, such conditions become the analytical tool to identify the boundary behaviour of the process: proving that every $f \in D(\mathcal{L})$ satisfies a given boundary condition is analytically equivalent to verifying that the process exhibits the corresponding probabilistic boundary behaviour (see, e.g., \cite[Secs.~2.5 and 3.5]{Mandl1968AnalyticalTO}).

In particular, the (reflecting) elastic 
behaviour at $0$ is identified 
through a \emph{Robin condition} on 
$D(\mathcal{L})$, which leads to the following definition.

\begin{definition}[Feller $(\mu, \sigma)$--reflected elastic diffusion]\label{def: Feller RED}
Let $\hat{X}$ be a Feller diffusion on $\mathcal{I} := [0,\infty)$ with infinitesimal parameters $(\mu, \sigma^2)$, scale function $\hat{S}$, and infinitesimal generator $\hat{\mathcal{L}}$. 
We say that $\hat{X}$ realizes the $(\mu, \sigma)$--\emph{Reflected Elastic Diffusion} (shortly, $(\mu, \sigma)$--RED) at $0$ with parameter $k \in (0,\infty)$ if every $f \in \mathcal{D}(\hat{\mathcal{L}})$ satisfies the \emph{Robin condition}
\begin{equation}
    \label{eq: robin_cond}
    \dfrac{df}{d\hat{S}}^+(0+) = k\, f(0). 
\end{equation}
\end{definition}

We now define the building-block for the construction of a RED. 
\begin{definition}[Weak solution of reflected SDE]\label{def: weak_sol_refl_sde}
Let $\mu:[0,\infty)\to\mathbb{R}$ and $\sigma:[0,\infty)\to\mathbb{R}$ be Borel measurable functions.
For $x\ge0$, a \emph{weak solution} of the reflected SDE
\begin{equation}\label{eq: refl_sde}
    \ud X_t = \mu(X_t)\,\ud t + \sigma(X_t)\,\ud W_t + \ud L_t^0, \qquad X_0 = x
\end{equation}
is a tuple $(\Omega, \mathcal{F}, (\mathcal{F}_t)_{t \ge 0}, \mathbb{P}, W, X, L^0)$, where
$(\Omega, \mathcal{F}, (\mathcal{F}_t)_{t\ge0}, \mathbb{P})$ is a filtered probability space with $\mathcal{F}$ $\mathbb{P}$-complete and $\mathcal{F}_0$ containing all $\mathbb{P}$-null sets of $\mathcal{F}$,
supporting a Brownian motion $W$ and continuous adapted processes $(X,L^0)$ such that
$\mathbb{P}$-a.s.
$$X_t \ge 0 \quad \forall\, t \ge 0, \qquad
L^0 \text{ is non-decreasing with } L_0^0 = 0, \qquad
\int_0^\infty \mathbf{1}_{\{X_t > 0\}}\,\ud L_t^0 = 0,$$
and \eqref{eq: refl_sde} holds $\mathbb{P}$-a.s.

We say that \eqref{eq: refl_sde} admits a \emph{unique weak solution} if a weak solution exists for every $x \ge 0$ and any two weak solutions starting from the same initial point share the same law on $C([0,\infty))$.
\end{definition}

Throughout the paper, we operate under the following set of assumptions. 
\begin{assumptions}\label{assumpt: refl_sde}
The coefficients $(\mu, \sigma)$ are continuous, and $\sigma >0$ on $[0,\infty)$.
For every initial state $x \ge 0$, the reflected SDE \eqref{eq: refl_sde} admits a unique weak solution
in the sense of Definition~\ref{def: weak_sol_refl_sde}, with associated strong Markov process
$X:=(X_t)_{t\ge0}$, laws $(\mathbb{P}_x)_{x\ge0}$,
and corresponding expectations $(\mathbb{E}_x)_{x\ge0}$.
We further assume that
\begin{enumerate}[(i)]
\item \label{assumpt: refl_sde_i} $\exists\, T>0,\ p>1$ such that
$$\mathbb{E}_0\left[\left(\int_0^T |\mu(X_u)|\, \ud u\right)^{2p}\right] < \infty, \quad
\mathbb{E}_0\left[\left(\int_0^T \sigma^2(X_u)\, \ud u\right)^{p}\right] < \infty.$$
\item \label{assumpt: refl_sde_ii} For every $t > 0$ and $M > 0$, $\lim_{x \to \infty} \mathbb{P}_x(X_t < M) = 0$.
\end{enumerate}
\end{assumptions}

The process $X$ takes values in $\mathcal{I} = [0,\infty)$, 
which we equip with the Borel $\sigma$-algebra $\mathcal{B}(\mathcal{I})$ in forming the state space $(\mathcal{I}, \mathcal{B}(\mathcal{I}))$ of $X$. In the following, $b\mathcal{B}(\mathcal{I})$ represents the set of bounded Borel functions on $\mathcal{I}$.
 We also denote by $(\mathcal{F}^X_t)_{t\ge0}$ the natural filtration of $X$, $\mathcal{F}^X_t:=\sigma(X_u,\,u\le t)$,
 augmented, for each $x\ge0$, with the $\mathbb{P}_x$-null sets.

A few comments regarding sufficient conditions for Assumption~\ref{assumpt: refl_sde}-\eqref{assumpt: refl_sde_i}--\eqref{assumpt: refl_sde_ii} are in order.

\begin{remark}[Sufficient conditions for Assumption~\ref{assumpt: refl_sde}]
     Both Assumptions~\ref{assumpt: refl_sde}-\eqref{assumpt: refl_sde_i} and \eqref{assumpt: refl_sde_ii} hold whenever $\mu$ and $\sigma$ satisfy a \emph{linear growth condition}, i.e., if there exists $K>0$ such that 
\begin{equation}\label{eq: linear_growth}
        |\mu(x)| + |\sigma(x)| \le K(1+|x|), \quad \forall x \ge 0.   
\end{equation}
    Indeed, under linear growth, Hölder's inequality reduces the integrability conditions in
    Assumption~\ref{assumpt: refl_sde}-\eqref{assumpt: refl_sde_i} to the problem of establishing
    $\mathbb E_0[\sup_{s\le T} X_s^{2p}]<\infty$. For $X_0=x=0$, writing $X=Y+L^0$ with
    $Y_t :=\int_0^t\mu(X_u)\,\ud u+\int_0^t\sigma(X_u)\,\ud W_u$, the Skorokhod formula (cf.~\cite[Ch.~VI, Lemma 2.1]{revuz2013continuous})
    $L^0_t=\sup_{s\le t}(-Y_s)^+$ gives 
    \begin{equation}\label{eq: X_vs_Y}
    \sup_{s\le t}X_s\le 2\sup_{s\le t}|Y_s|,\qquad \forall t\ge0.
\end{equation}
This allows to establish Assumptions~\ref{assumpt: refl_sde}-\eqref{assumpt: refl_sde_i} from \eqref{eq: linear_growth} via standard arguments, which are briefly reported here. 
Let $\tau_n:=\inf\{t\ge0:X_t\ge n\}$, so that $\tau_n\uparrow\infty$ a.s.\ since $X$ is a
continuous process on $[0,\infty)$. For any $t\le T$,  set
\[
    \phi_n(t):=\mathbb E_0\Big[\sup_{s\le t\wedge\tau_n}|Y_s|^{2p}\Big].
\]
Then
\[
    \phi_n(t)\le 2^{2p-1}\Big(
    \mathbb E_0\Big[\Big(\int_0^{t\wedge\tau_n}|\mu(X_u)|\,\ud u\Big)^{2p}\Big]
    +\mathbb E_0\Big[\sup_{s\le t\wedge\tau_n}\Big|\int_0^{s}\sigma(X_u)\,\ud W_u\Big|^{2p}\Big]\Big),
\]
and by
H\"older's inequality, linear growth and Burkholder–Davis–Gundy inequality (with constant $C_p$) we obtain
\begin{equation}\label{eq: phi_finiteness}
    \phi_n(t)\le C_1\,
    \mathbb E_0\Big[\int_0^{t\wedge\tau_n}(1+X_u)^{2p}\,\ud u\Big] \le C_1 2^{2p-1} \Big( T + \mathbb E_0\Big[\int_0^{t \wedge \tau_n} X_u^{2p}\, \ud u \Big]\Big),
\end{equation}
where $C_1: =  2^{2p-1}K^{2p}\big(T^{2p-1}+C_pT^{p-1}\big)$.
On $[0, \tau_n]$ we have $X \in [0,n]$, hence \eqref{eq: phi_finiteness} implies $\phi_n(t)<\infty, \forall t \le T$. Moreover, using $\mathbf 1_{\{u\le\tau_n\}}X_u\le\sup_{r\le u\wedge\tau_n}X_r$ and 
swapping integrals
$$
\phi_n(t) \le C_1 2^{2p-1}\Big(T+\int_0^t\mathbb E_0\Big[\sup_{r\le u\wedge\tau_n}X_r^{2p}\Big]\ud u\Big) \le C_1 2^{2p-1}\Big(T+2^{2p}\int_0^t\phi_n(u)\,\ud u\Big),
$$
where the last inequality follows from \eqref{eq: X_vs_Y}. 
Thus, for $t\le T$ we conclude
\[
    \phi_n(t)\le C\Big(1+\int_0^t\phi_n(u)\,\ud u\Big),
\]
with $C=C(p,K,T)$ independent of $n$.
Since $\phi_n\ge 0$ is non-decreasing and $\phi_n(T)<\infty$, it is a non-negative bounded Borel function on $[0,T]$,
so Gr\"onwall's lemma applies (cf.~\cite[Appendix, \S1]{revuz2013continuous}), yielding $\phi_n(t)\le Ce^{CT}, t\le T$.
Finally,  by \eqref{eq: X_vs_Y} and monotone convergence 
$$ 
    \mathbb E_0\Big[\sup_{s\le T}X_s^{2p}\Big]  \le 2^{2p}\mathbb E_0[\sup_{s\le T} |Y_s|^{2p}] =  2^{2p}\lim_{n\to \infty} \phi_n(T) \le 2^{2p}Ce^{CT}<\infty.
    $$

    Similarly, Assumption~\ref{assumpt: refl_sde}-\eqref{assumpt: refl_sde_ii} can be deduced from
linear growth. Let $g(y):=\frac{1}{1+y}$, so that $0<g\le1$ for $y\ge 0$, $g'(y)=-(1+y)^{-2}$,
$g''(y)=2(1+y)^{-3}$ and $g'(0)=-1$. By Itô's formula,
\[
    g(X_t) = g(x) + \int_0^t \bigg(g'(X_s)\mu(X_s) + g''(X_s) \frac{\sigma^2(X_s)}{2} \bigg) \,\ud s + \int_0^t g'(X_s)\sigma(X_s)\,\ud W_s
    - L_t^0 ,
\]
where we used that $L^0$ increases only when $X=0$, so that
$\int_0^t g'(X_s)\,\ud L^0_s=g'(0)L^0_t=-L^0_t$. By linear growth,
$|\mu(y)|(1+y)^{-2}\le K g(y)$ and $\sigma^2(y)(1+y)^{-3}\le K^2 g(y)$; moreover
$|g'(y)\sigma(y)|\le K$ for $y\ge 0$, so the stochastic integral is a true martingale. Dropping
$-L_t^0\le0$ and taking expectations, with $C:=K(K+1)$,
\[
    \mathbb{E}_x[g(X_t)] \le g(x) + C\int_0^t \mathbb{E}_x[g(X_s)]\,\ud s,
\]
and Grönwall's lemma yields $\mathbb{E}_x[g(X_t)]\le e^{Ct}/(1+x)$. Finally, by Markov's inequality,
\[
    \mathbb{P}_x(X_t<M)=\mathbb{P}_x\!\left(g(X_t)>\tfrac1{1+M}\right)
    \le (1+M)\,\mathbb{E}_x[g(X_t)] \le \frac{(1+M)e^{Ct}}{1+x}\xrightarrow{x\to\infty}0 .
\]
\end{remark}

\medskip 

For a fixed $x_0\in(0,\infty)$, the scale density $S'$ and speed density $m$ associated to the process $X$ are
defined by
\begin{equation}\label{eq: scale_speed_X}
    S'(x) = \exp\!\left\{-\int_{x_0}^x \frac{2\mu(z)}{\sigma^2(z)}\,\ud z\right\},  \qquad m(x) = \frac{2}{\sigma^2(x)\,S'(x)}, \qquad x \in (0,\infty).
\end{equation}
We write $S(x) = \int_{x_0}^x S'(z)\,\ud z$ for the \emph{scale function} of $X$ and
$m(\ud x) = m(x)\,\ud x$ for its \emph{speed measure}, for $x\in(0,\infty)$.

For future reference, we now
record some important properties of $S'$ and $m$, which in turn yield the \emph{regularity} of the process (cf.~\cite[Ch.~V, Def.~45.2]{rogers2000diffusions} or \cite[Ch.~2.2, Def.~9]{freedman2012brownian}).
\begin{remark}[]\label{rem: scale_speed}
Since $\mu,\sigma^2 \in C([0,\infty))$ with $\sigma>0$, 
the ratio $2\mu/\sigma^2$ is continuous and thus locally bounded on $[0,\infty)$. It follows that $S'(0+) := \lim_{x\downarrow 0}S'(x)$ exists in $(0,\infty)$, so that $S'$ extends continuously to $x=0$ by setting $S'(0) := S'(0+)$, yielding $0 < S' \in C([0,\infty))$. This implies $0 < m \in C([0,\infty))$ as well. Furthermore, by the properties of $S'$ we also deduce that $(X_t)_{t\ge0}$ is a \emph{regular} process on $[0,\infty)$ in the sense of \cite{rogers2000diffusions}, i.e.
\[
    \mathbb{P}_x(\tau_y < \infty) > 0 \qquad \text{for all } x \in (0,\infty), y \in [0, \infty) ,
\]
where $\tau_y := \inf\{t > 0 : X_t = y\}$ is the first hitting time of $y$. To see this, for $x,y>0$ note that $\mathbb{P}_x(\tau_y < \infty)$ can be bounded from below by restricting the process to a compact subinterval containing $x$:
$$\mathbb{P}_x(\tau_y < \infty) \ge \begin{cases} \mathbb{P}_x(\tau_y < \tau_a) = \dfrac{S(x) - S(a)}{S(y) - S(a)}, & \text{for any fixed } a \in (0, x) \text{ if } x \le y, \\ \mathbb{P}_x(\tau_y < \tau_b) = \dfrac{S(b) - S(x)}{S(b) - S(y)}, & \text{for any fixed } b \in (x, \infty) \text{ if } x > y, \end{cases}$$
(cf.~\cite[Ch.~2, No.~4]{Borodin}).
Since $S'>0$, 
both numerators and denominators above are strictly positive, implying $\mathbb{P}_x(\tau_y < \infty)>0$ for $ x,y\in(0,\infty)$. 

Finally, for $x>0=y$, let $b \in (x,\infty) $ and consider $\tau_{0+} := \lim_{y\downarrow 0} \tau_y$, the limit existing since $\tau_y$ is a non-increasing function of $y$. For every sample path $\omega$ starting at $x > 0$, $\tau_0(\omega) = \tau_{0+}(\omega)$ holds by continuity of trajectories:   $\tau_0 (\omega) \ge  \tau_{0+}(\omega)$ follows trivially by $\tau_0(\omega) \ge \tau_y(\omega), \forall y \in (0,x)$, yielding $\tau_0(\omega) = \tau_{0+}(\omega)$ if $ \tau_{0+}(\omega) = \infty$. Instead, if $ \tau_{0+}(\omega) < \infty$, then $X_{\tau_{0+}} = \lim_{y \downarrow 0} X_{\tau_y} = 0$, implying $\tau_{0+}(\omega) \ge \tau_0(\omega)$. Thus 
\begin{equation}\label{eq: tau_0_1}
    \mathbb{P}_x(\tau_0 <\infty) \ge \mathbb{P}_x(\tau_0 < \tau_b) = \mathbb{P}_x(\tau_{0+} < \tau_b), \quad 0<x<b.
\end{equation}
Since $\tau_{0+}=\tau_b<\infty$ is impossible by continuity of trajectories (i.e., $X_{\tau_{0+}}=0\neq b=X_{\tau_b}$), we have $\{\tau_{0+}=\tau_b\}=\{\tau_{0+}=\tau_b=\infty\}$, so \begin{equation} \mathbb{P}_x(\tau_{0+} < \tau_b) = \mathbb{P}_x(\tau_{0+} \le \tau_b) - \mathbb{P}_x(\tau_{0+} = \tau_b = \infty).
\end{equation}
By  ~\cite[Ch.~15, Eq.~(6.6)]{karlin1981second}, for every $y\in(0,x)$, $\mathbb{E}_x[\tau_y\wedge \tau_b] \le \int_0^b S'(z) \,\ud z \int_0^b m(z)\,\ud z=: C_b<\infty$ since $S',m \in C([0,\infty))$. This implies $\mathbb{P}_x(\tau_y = \tau_b = \infty) = 0$ for every $y\in(0,x)$, and monotone convergence yields $\mathbb{E}_x[\tau_{0+}\wedge \tau_b] = \lim_{y \downarrow 0}\mathbb{E}_x[\tau_y\wedge \tau_b] \le C_b<\infty$, hence $\mathbb{P}_x(\tau_{0+} = \tau_b = \infty) = 0$. Thus, for a sequence  $(y_n)_{n \in \mathbb{N}}\subset(0,x)$ such that $y_n \downarrow 0$ as $n \to \infty$, using $\tau_{0+}=\sup_n\tau_{y_n}$ and $\mathbb{P}_x(\tau_{y_n}=\tau_b)=0$,
\begin{align}\label{eq: tau_0_2}
 \mathbb{P}_x(\tau_{0+} < \tau_b)  &= \mathbb{P}_x(\tau_{0+} \le \tau_b) = \mathbb{P}_x\left( \bigcap_{n \in \mathbb{N}} \{\tau_{y_n} \le \tau_b\} \right) = \mathbb{P}_x\left( \bigcap_{n \in \mathbb{N}} \{\tau_{y_n} < \tau_b\} \right) \\&= \lim_{n  \to \infty} \mathbb{P}_x(\tau_{y_n} < \tau_b) = \lim_{n\to \infty}  \dfrac{S(b) - S(x)}{S(b) - S(y_n)} >0    
\end{align}
since $\lim_{n \to  \infty}\big(S(b) - S(y_n)\big) = \int_0^b S'(z) \,\ud z\in(0,\infty)$ by $0<S' \in C([0,\infty))$ (see also ~\cite[Ch.~15, Lemma $6.1$(i)]{karlin1981second} and the preceding discussion). 
 By \eqref{eq: tau_0_1}-\eqref{eq: tau_0_2} we  conclude $\mathbb{P}_x(\tau_0 <\infty) > 0$. The process $X$ is then regular.
\end{remark}

Under Assumption~\ref{assumpt: refl_sde}, consider a solution $(X, L^0)$ to \eqref{eq: refl_sde} and set $\lambda := k S'(0+)$ for some fixed $k \in (0,\infty)$. 
Given an exponential random variable $\mathrm{E} \sim \operatorname{Exp}(\lambda)$ independent of $\mathcal{F}_\infty := \sigma\Big(\bigcup_{t \ge 0} \mathcal{F}_t\Big)$, 
we define the random time $\zeta := \inf \{ t \ge 0 : L_t^0 \ge E \}$ and denote by $\widetilde X := (\widetilde X_t)_{t \ge 0}$ the reflected process subject to random-time killing at $\zeta$:
\[
\widetilde X_t: = 
\begin{cases}
    X_t & t<\zeta\\
    \dagger & t\ge \zeta
\end{cases}
\]
where $\{\dagger\}$ is an isolated cemetery state, so $\zeta$ represents the \emph{lifetime} of $\widetilde X$. 
By convention, we extend any function $f$ on $\mathcal{I}$ to $\mathcal{I} \cup \{\dagger\}$ by setting $f(\dagger) = 0$.

For $x\in(0,\infty)$, we denote by $\widetilde S(x)$ and $\widetilde m(\ud x)$ the scale function and speed measure of $\widetilde X$, with corresponding densities $\widetilde S'(x)$ and $\widetilde m(x)$.

\begin{remark}\label{rem: tildeX_scale_speed}

Since $\widetilde X$ is governed by the same infinitesimal coefficients $(\mu,\sigma^2)$ as $X$, its scale and speed densities coincide with those of $X$ in the interior $ (0,\infty)$ of the state space $\mathcal{I}$ (where densities are defined: see \eqref{eq: scale_speed_X}). More precisely, the choice of an arbitrary reference point $x_0 \in (0,\infty)$ introduces a degree of freedom, causing the scale and speed densities of $X$ and $\widetilde X$ to coincide only up to a positive multiplicative constant (and its reciprocal). In what follows, we fix a common reference point $x_0 \in (0,\infty)$, yielding $\widetilde S'(x)=S'(x)$ and $\widetilde m(x)=m(x)$, hence $\widetilde S(x) = S(x)$ and $\widetilde m(\ud x) = m(\ud x)$, for $x\in(0,\infty)$.
\end{remark}

We now state the main result: 
\begin{theorem}\label{thm: main}
If Assumption~\ref{assumpt: refl_sde} holds, the process $\widetilde X$ is the unique (in law) Feller $(\mu,\sigma)$-reflected elastic diffusion at $0$ with parameter $k$.
\end{theorem}

\section{Properties of the process $\widetilde X$}\label{sect: tildeX_properties}

This section collects the key 
properties of the process $\widetilde X$ needed to prove the main theorem. Specifically, we characterize the nature of the killing time $\zeta$ and show that the killed process $\widetilde X$ inherits the strong Markov property. Next, we leverage general diffusion theory to show that its transition semigroup possesses the Feller property. Finally, we deduce an explicit differential representation for the infinitesimal generator of $\widetilde X$. 

To start, we note that the random killing-time $\zeta$ introduced in Section \ref{sect: setting_main_result} represents a stopping time with respect to the natural filtration of $\widetilde X$, which we denote by $( \mathcal{F}^{\widetilde X}_t)_t$. Indeed 
$$ \{ \zeta\le t\} = \{\widetilde X_t \notin \mathcal{I} \} \in \mathcal{F}^{\widetilde X}_t, \quad \forall t\ge 0.$$

 Observe also that, denoting by $\bar F_{\mathrm{E}}$ the survival function of the random variable $\mathrm{E}$, the \emph{conditional survival probability} reads  
 \begin{equation}
 \label{eq: surv_prob}
     \mathbb{P}_x(t<\zeta |\mathcal{F}_t) = \mathbb{P}_x(L_t^0<E |\mathcal{F}_t) = \bar F_E(L_t^0) =  e^{-\lambda L_t^0} 
 \end{equation}
  by conditioning under independence (see, e.g., \cite[Ch.~4.2, Lemma 4.1]{baldi2017stochastic}), using that $L_t^0$ is $\mathcal{F}_t$-measurable (cf.~Def.~\ref{def: weak_sol_refl_sde}) and  $\mathrm{E} \independent
  \mathcal{F}_\infty$, where $\mathcal{F}_t \subseteq \mathcal{F}_\infty$.

To investigate the Markovianity of the process $\widetilde X$, it is of interest to introduce the notion of \emph{sub-process}: given a LCCB space $E$ and a couple of Markov processes  $X_1$ and $X_2$ on the same state space $(E, \mathcal{B}(E))$ with semigroups $(P^1_t)_{t\ge 0}$ and $(P^2_t)_{t\ge 0}$ , $X_2$ is a said to be a sub-process of $X_1$ if $P^2_t f(x) \le P_t^1 f(x), \forall x\in E, f \in b \mathcal{B}(E)$ such that (s.t.) $f\ge 0$ (we recall that any $f$ is taken to vanish on $\{\dagger\}$). 
In this respect, the process $\widetilde X$ can be qualified as the sub-process obtained via killing of the Markov process $X$ by a \emph{multiplicative functional} (MF) (we refer to \cite[Ch. III]{blumenthal2007markov} for a complete discussion on the topic); its Markovianity follows 
from this construction, as verified in the proof of the next lemma.

\begin{lemma}\label{lem: strong_markov}
 The process $\widetilde X$ is strong Markov with transition semigroup 
   \begin{equation}\label{eq: semigroup}
       P_t f(x) = \mathbb{E}_x[e^{-\lambda L_t^0} f(X_t)] = \mathbb{E}_x[f(\widetilde X_t)], \qquad  t\ge 0,\ f \in b \mathcal{B}(\mathcal{I}), \ x \in \mathcal{I}.
   \end{equation}
\end{lemma}
\begin{proof}
Recalling $\sigma(0)>0$ by Assumption~\ref{assumpt: refl_sde}, it is well known
(see, e.g., \cite[Ch.~1.3]{pilipenko2014reflection}) that $L^0$ almost surely coincides with half
of the semimartingale right local time of $X$ at $0$. That is, denoting
$A_t(\epsilon):=\int_0^t \mathbf{1}_{[0,\epsilon]}(X_u)\,\sigma^2(X_u)\,\ud u$, for every $x\ge0$,
\begin{equation}\label{eq: L_as_semimart_local_time}
    \mathbb{P}_x\Big( L_t^0 = \lim_{\epsilon\downarrow0}\frac{A_t(\epsilon)}{2\epsilon} \ \ \text{for all } t\ge0 \Big)=1 .
\end{equation}
 In the remainder of the proof, $\Omega=C([0,\infty);[0,\infty))$ denotes the
 canonical space and $\theta_t$ the shift operator,
 $(\theta_t\omega)(s):=\omega(t+s)$ for all $s\ge 0$.
 For a random time $\rho$, the shift $\theta_\rho$ is defined on $\{\rho<\infty\}$ 
 by
 $(\theta_\rho\omega)(s):=\omega(\rho(\omega)+s)$ for all $s\ge 0$
 (cf.~\cite[Ch.~1, Eq.~7.4]{blumenthal2007markov}).
 
 Let
$
G:=\big\{\omega\in\Omega:\ \lim_{\epsilon\downarrow0}\frac{A_u(\epsilon)(\omega)}{2\epsilon}\ \text{exists for all } u\ge0\big\}.
$
By \eqref{eq: L_as_semimart_local_time}, $\mathbb{P}_x(G)=1$ for every $x\ge0$.
Without loss of generality, we may and do replace $L^0$ with the version
\begin{equation}\label{eq: L_ind}
    L^0_t(\omega):=\begin{cases}
\displaystyle\lim_{\epsilon\downarrow0}\frac{A_t(\epsilon)(\omega)}{2\epsilon}, & \omega\in G,\\[2mm]
0, & \omega\notin G,
\end{cases}\qquad t\ge0,
\end{equation}
which, by \eqref{eq: L_as_semimart_local_time}, is $\mathbb{P}_x$-indistinguishable from the original one for every $x\ge0$. For every $\omega\in\Omega$, $L^0$ is non-decreasing with $L^0_0=0$,
and it is $\mathbb{P}_x$-a.s.\ continuous. Moreover, it is $(\mathcal{F}^X_t)_t$-adapted by completeness of the natural filtration
(cf.~Section~\ref{sect: setting_main_result}).
Finally, $L^0$ qualifies as an \emph{additive functional} of $X$
(cf.~\cite[Ch.~4, Def.~1.1]{blumenthal2007markov}) provided that, for any $s,t\ge0$,
\begin{equation}\label{eq: af}
    L^0_{t+s}=L^0_s+L^0_t\circ\theta_s
\end{equation}
holds $\mathbb{P}_x$-a.s.\ for every $x\ge0$.
Fix $\omega\in G$ and $s,t\ge0$. Since 
\begin{equation}\label{eq: A_1}
    A_{t+s}(\epsilon) - A_s(\epsilon) = \int_s^{t+s}\mathbf{1}_{[0,\epsilon]}(X_u)\,\sigma^2(X_u)\,\ud u =\int_0^{t}\mathbf{1}_{[0,\epsilon]}(X_{s+u})\,\sigma^2(X_{s+u})\,\ud u 
= A_t(\epsilon) \circ \theta_s 
\end{equation} 
holds pathwise
for every $\epsilon>0$, 
 we have
\begin{equation}\label{eq: A_2}
    \frac{A_t(\epsilon)(\theta_s\omega)}{2\epsilon}
=\frac{A_{t+s}(\epsilon)(\omega)-A_s(\epsilon)(\omega)}{2\epsilon}
\xrightarrow[\epsilon\downarrow0]{} L^0_{t+s}(\omega)-L^0_s(\omega),
\end{equation}
where the right-hand side of the equality converges because $\omega\in G$, implying that also the left-hand side converges as $\epsilon \downarrow 0$: since $t\ge0$ is arbitrary, we conclude
$\theta_t\omega\in G$, 
so that, by \eqref{eq: L_ind},
$(L^0_t\circ\theta_s)(\omega)=L^0_{t+s}(\omega)-L^0_s(\omega)$ holds on the $\mathbb{P}_x$-a.s.~event $G$, and  \eqref{eq: af} follows.

We now extend \eqref{eq: af} to any $(\mathcal{F}_t)_t$-stopping time $\tau<\infty$. Fix $\omega\in G$ and set $s:=\tau(\omega)$. Note that \eqref{eq: A_1} continues to hold with this choice of $s$. Moreover, since $\omega\in G$, the limit $\lim_{\epsilon\downarrow0}A_u(\epsilon)(\omega)/(2\epsilon)$ exists at $u=\tau(\omega)$ and at $u=\tau(\omega)+t$ for every $t\ge0$. It follows that for $s= \tau(\omega)$ the right-hand side limit in \eqref{eq: A_2} exists for every $t\ge0$. Hence $\theta_\tau\omega\in G$ and, by \eqref{eq: L_ind},
\begin{equation}\label{eq: strong_af}
    L^0_{t + \tau} = L^0_\tau + L^0_t \circ \theta_\tau
\end{equation}
holds $\mathbb{P}_x$-a.s.~for all $t\ge 0$ and for any finite $(\mathcal{F}_t)_t$-stopping time $\tau$.

By virtue of the properties of $ L^0$, the process
$$M_t := e^{-\lambda  L_t^0}, \quad t\in [0,\infty), \quad M_\infty:= 0$$
defines a continuous (a.s.)~\emph{multiplicative functional} (MF) of $X$ (cf.~\cite[Ch.~3, Def.~1.1]{blumenthal2007markov}): $M_t$ is $\mathcal{F}^X_t-$measurable for all $t\ge 0$; $M_t(\omega) \in [0,1]$ for all $(t, \omega)$; and for any $t,s\ge 0$, $M_{t+s} = M_t (M_s \circ \theta_t), \mathbb{P}_x$ a.s. for every $x\ge 0$. Furthermore, $M$ qualifies as a \emph{strong MF} provided that 
\begin{equation}\label{eq: strong_MF}
\mathbb{E}_x[M_{t+\tau}f(X_{t+\tau})] = \mathbb{E}_x\big[M_{\tau} \mathbb{E}_{X_\tau}[M_{t}f(X_{t})] \big]
\end{equation}
holds for every $t \ge 0$, $x \in \mathcal{I}$, $f \in b\mathcal{B}(\mathcal{I})$, and every $(\mathcal{F}_t)_{t}$-stopping time $\tau$ (cf.~\cite[Ch.~4, Def.~3.11]{blumenthal2007markov}). If $\tau = \infty$ this is immediate by the convention $M_\infty = 0$. If $\tau<\infty$, by \eqref{eq: strong_af} 
$$\mathbb{E}_x\bigl[e^{-\lambda  L_{t+\tau}^0} f(X_{t+\tau})\bigr] = \mathbb{E}_x\bigl[e^{-\lambda  L_\tau^0 } \mathbb{E}_x\bigl[\big(e^{-\lambda  L_t^0}f(X_t)\big)\circ \theta_\tau |\mathcal{F_\tau}\bigr]\bigr] = \mathbb{E}_x\bigl[e^{-\lambda  L_\tau^0 } \mathbb{E}_{X_{\tau}}\bigl[\big(e^{-\lambda  L_t^0}f(X_t)\bigr]\bigr]$$
where the first equality holds by the tower property and the fact that $L^0$ is $(\mathcal{F}_t)_{t\ge0}$-progressively measurable
\footnote{
For each $\epsilon>0$, $(A_t(\epsilon))_{t}$ is adapted and continuous, hence
progressive; therefore $\limsup_{q\downarrow0,\,q\in\mathbb{Q}}A_t(q)/(2q)$ is
progressive. Since on $G$ the limit in \eqref{eq: L_ind} exists, it coincides with
this limsup, so that
$L^0_t=\mathbf{1}_G\,\limsup_{q\downarrow0,\,q\in\mathbb{Q}}A_t(q)/(2q)$
(with $0\cdot\infty=0$). 
Finally, $G^c$ is $\mathbb{P}_x$-null for every $x\ge0$,
so $\mathbf{1}_G$ is $\mathcal{F}_0\subseteq\mathcal{F}_t$-measurable for all $t$ by completeness of the reference filtration (cf.~Def.~\ref{def: weak_sol_refl_sde});  
being time-constant, $\mathbf{1}_G$ is also  progressive, and so is the product above.
}.
The second equality follows from
\cite[Ch.~1, Prop.~8.4]{blumenthal2007markov}, that is,  strong Markovianity applied to the bounded, $\mathcal{F}_t^X$-measurable random variable $e^{-\lambda  L_t^0}f(X_t)$.

The Markovianity of the process $\widetilde X$ then follows from standard results about sub-processes obtained by killing of Markov processes through MFs: in particular, following the \emph{canonical construction} in~\cite[Ch.~3, pp.~105--106]{blumenthal2007markov}, $\widetilde X$ is the sub-process of $X$ obtained via killing through $M$.
We now translate the notation and measures defined in \cite[pp.~105--106]{blumenthal2007markov} to the present setting, in order to apply their results to $\widetilde X$:
for each $\omega \in \Omega$, the  measure $\alpha_\omega$ on $[0,\infty]$ is given by $\alpha_\omega((t, \infty]) := M_t(\omega) = e^{-\lambda  L_t^0(\omega)}$ for $t\in[0,\infty]$  and  $\alpha_\omega(\{0\}) = 0$. For $x \in \mathcal{I}$,  the associated probability measure $\widetilde P_x$ on $\Omega \times [0, \infty]$ (the second coordinate carrying the killing time) is defined on product sets as
\[
    \widetilde P_x( A \times (t, \infty]) = \mathbb{E}_x\bigl[e^{-\lambda  L_t^0} \mathbf{1}_A \bigr], \qquad \text{for } A \in \mathcal{F}, \, t\ge 0,
\]
while for $x=\dagger$,  $\widetilde P_\dagger$ has unit mass on the event that the killing time equals $0$: under $\widetilde P_\dagger$ the process is killed instantaneously, whatever the underlying path of $X$, so that $\widetilde X_t=\dagger, \forall t\ge0$.

    With these measure choices, we can apply~\cite[Ch.~3, Thm.~3.3]{blumenthal2007markov} to conclude $\widetilde X$ is Markov with transition semigroup
    $$P_t f(x) = \mathbb{E}_x\bigl[e^{-\lambda L_t^0} f(X_t)\bigr], \qquad t\ge 0,\ x \in \mathcal{I},\ f \in b\mathcal{B}(\mathcal{I}),$$
    which, recalling \eqref{eq: surv_prob}, by the tower property can be equivalently written as    
    \begin{equation}
        P_tf(x) = \mathbb{E}_x\bigl[\mathbf{1}_{\{t<\zeta\}}\, f(X_t)\bigr] = \mathbb{E}_x\bigl[f(\widetilde X_t)\bigr], \qquad t\ge 0,\ x \in \mathcal{I},\ f \in b\mathcal{B}(\mathcal{I}).
    \end{equation}
    
Finally, since $M$ is a strong MF, the process $\widetilde X$ inherits the strong Markov property from $X$  (see \cite[Ch.~3, Prop.~3.12]{blumenthal2007markov}).

\end{proof}

Since $t \mapsto \widetilde{X}_t$ is continuous on $[0, \zeta)$, Lemma~\ref{lem: strong_markov}
implies that $(\widetilde{X}_t, t \in [0,\zeta))$ is a time-homogeneous, one-dimensional \emph{diffusion} on $\mathcal{I}$
in the sense of~\cite[Ch.\,2, No.\,1]{Borodin}.
In particular, since $\widetilde X = X$ up to killing and $X$ is regular (see Remark \ref{rem: scale_speed}), the process 
$\widetilde X$ inherits this regularity, as detailed in the remark below.

\begin{remark}[Regularity of $\widetilde X$]
Since $\widetilde X_t = \dagger$ for all $t \ge \zeta$, the regularity of $\widetilde X$ is established by assessing the strict positivity of
$$ \mathbb{P}_x(\tilde \tau_y < \infty) = \mathbb{P}_x(\tilde \tau_y < \zeta), \quad x \in (0,\infty),\ y \in [0,\infty), $$
where $\tilde \tau_y := \inf\{t > 0 : \widetilde X_t = y\}$ (recall $\tau_y = \inf\{t > 0 : X_t = y\}$ from Remark \ref{rem: scale_speed}). Since $\widetilde X = X$ on $[0,\zeta),$ pathwise,
\begin{equation}\label{eq:path_identity}
    \{\tilde \tau_y < \infty\} = \{\tau_y < \zeta\}.
\end{equation}
Moreover, $L^0_t = 0$ for all $t \le \tau_0$ and $\mathrm{E} > 0$, $\mathbb{P}_x\text{-a.s.}$, thus $\zeta \ge \tau_0$, $\mathbb{P}_x\text{-a.s.}$ In particular, on the event $\{\tau_0 < \infty\}$ we have $\zeta > \tau_0$, $\mathbb{P}_x\text{-a.s.}$ This holds trivially when $\zeta = \infty$; if $\zeta < \infty$, note that $L^0_{\zeta} = \mathrm{E}$ while $L^0_{\tau_0} = 0$, so $\tau_0 = \zeta$ contradicts $\mathrm{E} > 0$, which holds $\mathbb{P}_x\text{-a.s.}$ 

Then, we distinguish two cases:
\begin{itemize}
    \item[$(\mathrm{i})$] $0 < x \le y$: For $a \in (0,x)$, we have $\{\tau_y < \tau_a\} \subseteq \{\tau_y < \tau_0\}$ pathwise, and $\{\tau_y < \tau_0\} \subseteq \{\tau_y < \zeta\} = \{\tilde \tau_y < \infty\}$, $\mathbb{P}_x\text{-a.s.}$ Hence, $\mathbb{P}_x(\tilde \tau_y < \infty) \ge \mathbb{P}_x(\tau_y < \tau_a) > 0$ by the arguments in Remark~\ref{rem: scale_speed}.
    
    \item[$(\mathrm{ii})$] $0 \le y < x$: If $\tau_0 < \infty$, then $\zeta > \tau_0$ and $X_t = \widetilde X_t$ for all $t \le \tau_0$, $\mathbb{P}_x\text{-a.s.}$ Since $\tau_y \le \tau_0$, it follows that $\widetilde \tau_y = \tau_y$, $\mathbb{P}_x\text{-a.s.}$ If $\tau_0 = \infty$, then $\zeta = \infty$ and $X_t = \widetilde{X}_t$ for all $t \ge 0$, $\mathbb{P}_x\text{-a.s.}$, implying again that $\widetilde \tau_y = \tau_y$, $\mathbb{P}_x\text{-a.s.}$ We conclude that $\mathbb{P}_x(\tilde \tau_y < \infty) = \mathbb{P}_x(\tau_y < \infty) > 0$ by the regularity of $X$.
\end{itemize}
Hence, $\widetilde X$ is a regular process on $\mathcal{I}$.
\end{remark}

By the standard theory on one-dimensional regular diffusions, 
it is known that $\widetilde X$ admits a transition density with respect to its speed measure $\widetilde m(\ud y)$, which we denote by $\widetilde p$ (cf.~\cite[p.~149]{ito2012diffusion}). Thus, the semigroup \eqref{eq: semigroup} of $\widetilde X$ takes the integral form
\begin{equation}\label{eq: semigroup_m}
    P_tf(x)  
    =\int_\mathcal{I} f(y) \widetilde p(t,x,y) \widetilde m(\ud y) \qquad   t\ge 0,\ x \in \mathcal{I},\ f \in b\mathcal{B}(\mathcal{I}), 
\end{equation}
where $\widetilde p: (0,\infty) \times \mathcal{I} \times \mathcal{I} \to [0,\infty)$  can be taken to be jointly continuous in all variables for $t>0$ (see \cite[Ch.~V, Thm.~50.11]{rogers2000diffusions}).

We are now in a position to prove that $(P_t)_t$ is a strongly continuous, contraction semigroup of positive linear operators mapping $C_0(\mathcal{I})$ into itself, that is, it qualifies as a Feller semigroup  in the sense of \cite{revuz2013continuous}.

\begin{lemma}\label{lem: feller-dynkin}
The family $(P_t)_{t\ge 0}$ is a  Feller semigroup.
\end{lemma}
\begin{proof}
For each fixed $t\ge0$, $P_t$ is a positive linear operator by the properties of expectation, and $P_0=I$ holds trivially since $L^0_0=0$.

The semigroup property of $(P_t)_{t\ge0}$, already stated in Lemma~\ref{lem: strong_markov}, is a special case of the general fact that the transition operator generated by a multiplicative functional of a Markov process is a semigroup (see \cite[Ch.~3, Sect.~1]{blumenthal2007markov}); in our setting it can also be verified directly, as follows: for any $s,t\ge0$, $x\in\mathcal I$, $f\in b\mathcal B(\mathcal I)$, taking $L^0$ as in \eqref{eq: L_ind} we have
\begin{align*}
 P_{t+s}f(x) &= \mathbb E_x\big[e^{-\lambda L^0_{t+s}}f(X_{t+s})\big] 
=\mathbb E_x\big[e^{-\lambda  L^0_s}\,\mathbb E_x\big[(e^{- \lambda  L^0_t}f(X_t))\circ\theta_s\,\big|\,\mathcal F_s\big]\big] \\
  &= \mathbb E_{x}\big[e^{-\lambda  L^0_s}\,\mathbb E_{X_s}\big[e^{- \lambda  L^0_t}f(X_t)\big]\big] = P_s P_tf(x)
\end{align*}
where the second equality follows from the additivity
$ L^0_{t+s}= L^0_s+ L^0_t \circ\theta_s,$ $\mathbb P_x$ a.s., the tower property,
and the $\mathcal F_s$-measurability of $e^{-\lambda L^0_s}$;
the third equality is an application of the Markov property of $X$ 
(in the sense of \cite[Ch.~6.1, Prop.~6.1]{baldi2017stochastic}, time-homogeneous form) to
$(e^{-\lambda  L^0_t}f(X_t))\circ\theta_s=e^{-\lambda( L^0_{t+s}- L^0_s)}f(X_{t+s})$,
which is bounded and $\sigma(X_u,\,u\ge s)$-measurable (as usual, up to completion with $\mathbb{P}_x$-null sets, $x\ge0)$ by the construction in \eqref{eq: L_ind}.

For any $f \in C_0(\mathcal{I})$ and $x\in\mathcal{I}, |P_tf(x)| \le \mathbb{E}_x\left[|f(X_t) \right|] \le ||f||_\infty$, thus $||P_tf||_\infty\le  ||f||_\infty$, proving that $P_t$ is a contraction operator. 

To show that $P_t$ maps $C_0(\mathcal{I})$ into $C_0(\mathcal{I})$, we first establish that $x \mapsto P_tf(x)$ is continuous for all $t \ge 0$. When $t=0$,  $P_0 = I$ implies that the continuity of the semigroup reduces to the one of $f$. For $t>0$, recall that \eqref{eq: semigroup_m} holds and consider a non-negative sequence $(x_n)_{n \in \mathbb{N}}$ with $x_n \to x$ as $n \to \infty$. 
Since $f \in C_0(\mathcal{I})$,  $\forall \epsilon >0, \exists C_\epsilon>0$ such that $\sup_{x\ge C_\epsilon} |f(x)| \le \epsilon/2$, thus, for $K_\epsilon: = [0, C_\epsilon]$
\begin{align*}
    |P_tf(x_n) - P_tf(x)| &= \bigg|\int_{K_\epsilon} f(y) \big(\widetilde p(t,x_n, y) - \widetilde p(t,x,y)\big)\widetilde m(\ud y)\\ 
    &  + \int_{\mathcal{I} \setminus K_\epsilon} f(y) \big(\widetilde p(t,x_n, y) - \widetilde p(t,x,y)\big)\widetilde m(\ud y) \bigg| \\
    &\le ||f||_\infty \int_{K_\epsilon}  \big|\widetilde p(t,x_n, y)- \widetilde p(t,x,y)\big|\widetilde m(\ud y) \\ 
    & +\frac{\epsilon}{2} \int_{\mathcal{I} \setminus K_\epsilon}  \big|\widetilde p(t,x_n, y) - \widetilde p(t,x,y)\big|\widetilde m(\ud y).
\end{align*}
By triangle inequality and non-negativity of the transition density,
\begin{align*}
    \int_{\mathcal{I} \setminus K_\epsilon}  \big|\widetilde p(t,x_n, y) - \widetilde p(t,x,y)\big|\widetilde m(\ud y) 
    \le \int_{\mathcal{I}}  \big(\widetilde p(t,x_n, y) + \widetilde p(t,x,y)\big)\widetilde m(\ud y) \le 2.
\end{align*}
Furthermore, the joint continuity of $(x,y) \mapsto \widetilde p(t,x, y)$ for any fixed $t>0$ implies that, on $K_\epsilon$, 
$$ \big|\widetilde p(t,x_n, y)- \widetilde p(t,x,y)\big| \le \widetilde p(t,x_n, y)+ \widetilde p(t,x,y) \le 2 \sup_{x'\in [0,N], y \in K_\epsilon} \widetilde p(t, x', y) = : 2\bar C_{N, \epsilon}< \infty$$
where $N > 0$ is chosen such that $\{x_n: n \in \mathbb{N}\} \cup \{x\} \subset [0, N]$. 
Note that $\int_{K_\epsilon}\bar C_{N, \epsilon} \widetilde m(\ud y) 
=\bar C_{N, \epsilon} \int_{K_\epsilon}  m(\ud y) < \infty$, 
 where the equality follows from Remark~\ref{rem: tildeX_scale_speed} and finiteness holds by the continuity of the speed density $m(y)$ on $[0,\infty)$ (see Remark~\ref{rem: scale_speed}).
 Hence, we apply the reverse Fatou's lemma to conclude 
$$ 0\le \limsup_{n \to \infty}   |P_tf(x_n) - P_tf(x)| \le  ||f||_\infty \int_{K_\epsilon}  \limsup_{n \to \infty} \big|\widetilde p(t,x_n, y)- \widetilde p(t,x,y)\big|\widetilde m(dy) +
 \epsilon = \epsilon,  $$
 and, since $\epsilon$ was arbitrary, 
 by letting it vanish we get 
 $$0 \le \liminf_{n\to \infty} |P_tf(x_n) - P_tf(x)| \le \limsup_{n \to \infty}   |P_tf(x_n) - P_tf(x)|  = 0 \Rightarrow \lim_{n \to \infty }P_t f(x_n) = P_tf(x).$$

Next, we verify that $P_t f$ vanishes at infinity: again, by $f \in C_0(\mathcal{I})$ 
    \begin{align*}
        |P_tf(x)| & \le \mathbb{E}_x \left [ |f(X_t)|\right]   =  \mathbb{E}_x[|f(X_t)|\mathbf{1}_{\{X_t< C_{\epsilon}\}}] +\mathbb{E}_x[|f(X_t)|\mathbf{1}_{\{X_t\ge C_{\epsilon}\}}] \\ & \le ||f||_\infty \mathbb{P}_x(X_t < C_\epsilon) + \frac{\epsilon}{2}.
    \end{align*}
By Assumption \ref{assumpt: refl_sde}-\eqref{assumpt: refl_sde_ii} it holds $\mathbb{P}_x(X_t < C_\epsilon) \rightarrow 0$ as $x\to \infty$. So, $\limsup_{x\to \infty} |P_t f(x)| \le \frac{\epsilon}{2}$. 
By arbitrariness of $\epsilon$ we conclude $\lim_{x \rightarrow \infty}P_tf(x) = 0$, confirming $P_tf \in C_0(\mathcal{I}), \forall f\in C_0(\mathcal{I})$. 

Finally, given the properties established above, the strong continuity of the semigroup, i.e., 
$$ \lim_{t\downarrow 0} || P_t f - f||_\infty = \lim_{t \downarrow 0} \sup_{x\in [0,\infty) }|P_tf(x) - f(x)| = 0,$$
follows from pointwise continuity (see, e.g., \cite[Ch. I, No. 11]{Borodin} or \cite[Ch. 3, Sect. 2, Prop. 2.4]{revuz2013continuous}), i.e., 
$$ \lim_{t \downarrow 0} P_t f(x) =  \lim_{t \downarrow 0} \mathbb{E}_x[e^{-\lambda L_t^0} f(X_t)] = f(x), \quad x\ge 0,$$
which holds by dominated convergence and continuity of the trajectories of $(X, L^0)$. \
\end{proof}

\begin{remark}
The notion of Feller semigroup is not uniform across the literature: what \cite{revuz2013continuous} calls a Feller semigroup corresponds to the \emph{Feller--Dynkin} semigroup of~\cite{Borodin}, where the term Feller semigroup is instead reserved for semigroups mapping $C_b(\mathcal{I})$ into itself. Throughout the paper, we stick to the definitions of \cite{revuz2013continuous}.
\end{remark}

We conclude that $\widetilde X$ is a \emph{Feller diffusion}. 
It follows that its \emph{infinitesimal generator}, that we denote by $(\widetilde{\mathcal{L}}, D(\widetilde{\mathcal{L}}))$, is defined by \eqref{eq: inf_gen_target}-\eqref{eq: dom_inf_gen_target} with $E = \mathcal{I}$. 

By Remark~\ref{rem: tildeX_scale_speed}, the scale and speed densities of $\widetilde X$ coincide with those of $X$ on $(0,\infty)$. Thus, $\widetilde X$ admits a speed density that is both strictly positive and continuous on $(0,\infty)$ (see Remark~\ref{rem: scale_speed}). Moreover, $\widetilde X$ is regular, hence it features no trap/absorbing points in $(0,\infty)$: for this class of diffusions, the approach developed in \cite[pp. 303-305]{karlin1981second} (based on the so-called \emph{Dynkin representation} of the generator) allows to show that, on $(0,\infty)$, the operator $\mathcal{\widetilde L}$ is the second-order differential operator
\begin{equation}
\label{eq: inf_gen_as_diff_op}
    \widetilde{\mathcal{L}} f(x) = \frac{\sigma^2(x)}{2} f''(x) + \mu(x) f'(x), \qquad f \in D(\widetilde{\mathcal{L}}), x \in (0,\infty).
\end{equation}

\section{Local time asymptotics}\label{sect: local_time_asymptotics}
The proof of the main theorem requires one last ingredient: the limiting behaviour, as $t\downarrow0$, of $L^0_t$, the local time of the reflected process $X$ at $\{0\}$. In this section we extend the findings of \cite[Lemma 15]{peskir2019continuity} and \cite[Lemma 4.1]{de2022stopping} in two directions. First, whereas \cite{peskir2019continuity, de2022stopping} provide a liminf lower bound on the expectation of $L^0$ for an \emph{unreflected} diffusion --- specifically, for the first coordinate of a two-dimensional elliptic diffusion in \cite{peskir2019continuity}, and for a one-dimensional diffusion in natural scale in \cite{de2022stopping} --- we establish the corresponding lower bound for a one-dimensional \emph{reflected} diffusion with drift $\mu$. Second, we provide a $\limsup$ upper bound on $\mathbb E[(L^0_t)^2]$, a second-moment estimate addressed in neither of the two papers.

\begin{proposition}\label{prop: local_time_asymptotics}
    Fix $\epsilon>0$, let ${K_\epsilon} := [0,\epsilon]$ and define $$0<\underline{\sigma}_{K_\epsilon} : = \inf_{z \in {K_\epsilon}} \sigma(z) \le    \sup_{z \in {K_\epsilon}} \sigma(z) =: \bar\sigma_{K_\epsilon}< \infty.$$ It holds 
    \begin{equation}\label{eq: liminf_bound}
        \liminf_{t \downarrow 0} t^{-1/2} \mathbb{E}_0[L_t^{0}] \ge \frac{\underline{\sigma}_{K_\epsilon}}{\sqrt{2\pi}}
    \end{equation}  
    and 
    \begin{equation}\label{eq: limsup_bound} \limsup_{t\downarrow 0} t^{-1} \mathbb{E}_0[(L_t^{0})^2] \le  2\bar \sigma_{K_\epsilon}^2.
    \end{equation}  
\end{proposition}
\begin{proof} 
    Fix $x= X_0 = 0$. By  Skorokhod’s lemma (cf.~\cite[Ch.~VI, Lemma 2.1]{revuz2013continuous}) 
    $$L_t^0 = \sup_{s \in [0,t]} \Big(-\int_0^s \mu(X_u) \ud u -  M_s\Big), \quad M_s := \int_0^s \sigma(X_u) \ud W_u$$
    with $(X_t)_{t\ge 0}$ evolving according to \eqref{eq: refl_sde}.
    For fixed $\epsilon>0$, let $\tau_\epsilon: = \inf\{t\ge 0: X_t \ge \epsilon \}$. By Assumption ~\ref{assumpt: refl_sde}, the function $\sigma$ is continuous and hence locally bounded, so $(M_{t \wedge \tau_\epsilon})_t$ is a continuous martingale. By similar arguments to those in \cite{peskir2019continuity, de2022stopping}, the Dubins-Schwarz theorem (cf.~\cite[Ch.~IV, Thm.~34.11]{rogers2000diffusions} with the trivial stopping sequence $\tau_n = n, n\in \mathbb{N}$) implies that the stopped version of $M$ can be represented as a time-changed Brownian motion, that is, $M_{t \wedge \tau_{\epsilon}} = \widetilde W_{\langle M\rangle_{t \wedge \tau_\epsilon}}, t\ge0$ where $\widetilde W$ is another Brownian motion and $\langle M\rangle _u = \int_0^u \sigma^2(X_s) ds$ is the quadratic variation of $M$. Since  $M_s = M_{s \wedge \tau_\epsilon}$ for $s \le t \wedge \tau_\epsilon$, it follows that
    \begin{equation}\label{eq: L_t_min_tau_eps}
         L^0_{t \wedge \tau_\epsilon}=   \sup_{s\le t \wedge \tau_\epsilon}\Big(- \int_0^s \mu(X_u)\ud u - \widetilde W_{{\langle M \rangle}_s}\Big).
    \end{equation} 
    
To prove \eqref{eq: liminf_bound}, note that $\mu$ is bounded by $c_\epsilon := \sup_{x \in [0,\epsilon]} \vert{}\mu(x)\vert{} < \infty$ on $[0,\epsilon]$ by continuity (cf.~Assumption \ref{assumpt: refl_sde}). Hence, using that $\langle M\rangle$ maps $[0, t \wedge \tau_\epsilon]$ continuously and monotonically onto $[0, \langle M\rangle_{t \wedge \tau_\epsilon}]$, we obtain
\begin{align}\label{eq: L_as_tilde_W}
    L^0_{t \wedge \tau_\epsilon} &\ge \sup_{s\le t \wedge \tau_\epsilon}\Big(- \int_0^s |\mu(X_u)|\ud u - \widetilde W_{{\langle M \rangle}_s}\Big) \ge   -\int_0^{t \wedge \tau_\epsilon} |\mu(X_u)|\ud u + \sup_{s\le t \wedge \tau_\epsilon} (-\widetilde W_{{\langle M \rangle}_s}) \\&\ge -c_\epsilon t +\sup_{s\le t \wedge \tau_\epsilon} (-\widetilde W_{{\langle M \rangle}_s}) = -c_\epsilon t +\sup_{u \le \langle M\rangle_{t \wedge \tau_\epsilon}} (-\widetilde W_{u}). 
\end{align}
Note that
     \begin{equation}\label{eq: M_bounds}
        \langle M\rangle_{t \wedge \tau_\epsilon} \in [\underline{\sigma}_{K_\epsilon}^2(t \wedge \tau_\epsilon), \bar \sigma^2_{K_\epsilon} (t \wedge \tau_\epsilon)].
    \end{equation} 
    Thus, by \eqref{eq: L_as_tilde_W}, \eqref{eq: M_bounds} and pathwise monotonicity of $L^0$, we obtain
    \begin{align}\label{eq: L_lower}
\mathbb{E}_0[L_t^0] \ge &  \mathbb{E}_0\left[L_{t \wedge \tau_\epsilon}^0\right] \ge \mathbb{E}_0\bigg[\sup_{u \le \underline{\sigma}_{K_\epsilon}^2(t \wedge \tau_\epsilon) }(-\widetilde W_u)\bigg]  - c_\epsilon t \ge   \mathbb{E}_0\bigg[\mathbf{1}_{\{t< \tau_\epsilon\}}\sup_{u \le \underline{\sigma}_{K_\epsilon}^2t }(-\widetilde W_u)\bigg]  -c_\epsilon t  \\=&  \mathbb{E}_0\bigg[\sup_{u \le \underline{\sigma}_{K_\epsilon}^2t }(-\widetilde W_u)\bigg] - \mathbb{E}_0\bigg[\mathbf{1}_{\{t\ge\tau_\epsilon\}}\sup_{u \le \underline{\sigma}_{K_\epsilon}^2t }(-\widetilde W_u)\bigg]  -c_\epsilon t \\
        \ge &   \mathbb{E}_0 \bigg[ \sup_{u \le \underline{\sigma}_{K_\epsilon}^2t } (-\widetilde W_u) \bigg] - \sqrt{ \mathbb{E}_0 \bigg[ \Big(\sup_{u \le \underline{\sigma}_{K_\epsilon}^2t } (-\widetilde W_u) \Big)^2 \bigg] } \sqrt{\mathbb{P}_0(t\ge\tau_\epsilon)} -c_\epsilon t,
    \end{align} 
    where the last bound follows by the Cauchy-Schwarz inequality. 
    By the symmetry and the reflection principle of Brownian motion,
    \begin{equation}\label{eq: bm_properties}
        \sup_{u \le \underline{\sigma}_{K_\epsilon}^2t }(-\widetilde W_u) \overset{d}{=} \sup_{u \le \underline{\sigma}_{K_\epsilon}^2t }(\widetilde W_u) \overset{d}{=}| \widetilde W_{\underline{\sigma}_{K_\epsilon}^2t}| \overset{d}{=} \underline{\sigma}_{K_\epsilon}\sqrt{t} |\widetilde W_1|, 
    \end{equation}
   where $\mathbb{E}[ |\widetilde W_1|] = \sqrt{\frac{2}{\pi}}$. Simple algebra stemming from \eqref{eq: L_lower} then yields
    $$   \mathbb{E}_0[L_t^0] \ge   \underline{\sigma}_{K_\epsilon}\sqrt{t} \bigg( \sqrt{\frac{2}{\pi}}- \sqrt{\mathbb{P}_0(t\ge\tau_\epsilon)} \bigg) -c_\epsilon t.$$
    Since $\mathbb{P}_0(\tau_\epsilon>0) = 1, \exists t_0>0$ s.t. $\forall t \in (0,t_0),\mathbb{P}_0(t\ge\tau_\epsilon) \le \frac{1}{2\pi}$. Thus, for all $t \in (0,t_0)$
     $$   \mathbb{E}_0[L_t^0] \ge   \frac{\underline{\sigma}_{K_\epsilon}}{\sqrt{2\pi}}\sqrt{t}-c_\epsilon t $$
    and \eqref{eq: liminf_bound} follows upon taking the $\liminf$ as $t\downarrow0$.
    
    It remains to prove \eqref{eq: limsup_bound}. Splitting the second-moment of $L^0$ based on $\tau_\epsilon$, exploiting the path monotonicity of $L^0$ and finally H\"older inequality with some conjugate exponents $p, q \in(1,\infty)$
    , we obtain
   \begin{align}\label{eq: L_upper}
       \mathbb{E}_0[(L_t^0)^2] = &  \mathbb{E}_0[\mathbf{1}_{\{t < \tau_\epsilon\}} (L_{t\wedge \tau_\epsilon}^0)^2] + \mathbb{E}_0[\mathbf{1}_{\{t \ge \tau_\epsilon\}} (L_{t}^0)^2] \le  \mathbb{E}_0[ (L_{t\wedge \tau_\epsilon}^0)^2] +\mathbb{E}_0[\mathbf{1}_{\{t \ge\tau_\epsilon\}}(L_{T}^0)^2] \\
       \le & \mathbb{E}_0[(L_{t\wedge \tau_\epsilon}^0)^2] +  \mathbb{E}_0[(L_{T}^0)^{2p}]^{\frac{1}{p}}\big(\mathbb{P}_0(t\ge\tau_\epsilon)\big)^{\frac{1}{q}}
   \end{align}
   for any fixed $T>t$. We now address each term in the last line separately. By arguments similar to those leading to \eqref{eq: L_as_tilde_W},
    $$ L_{t\wedge \tau_\epsilon}^0\le c_\epsilon  t + \sup_{u\le \langle M \rangle_{t\wedge \tau_\epsilon}} (- \widetilde W_u) \le c_\epsilon t + \sup_{u\le \bar \sigma_{K_\epsilon}^2 t} (-\widetilde W_u), $$
   where the last bound follows from \eqref{eq: M_bounds}. Hence, by $(a+b)^2\le 2a^2 + 2b^2$ and \eqref{eq: bm_properties}, 
   $$ \mathbb{E}_0[(L_{t\wedge \tau_\epsilon}^0)^2] \le 2(c_\epsilon t)^2 + 2\mathbb{E}_0\bigg[ \Big(\sup_{u\le \bar \sigma_{K_\epsilon}^2t} (-\widetilde W_u)\Big)^2\bigg] = 2(c_\epsilon t)^2 + 2\mathbb{E}_0\big[|\widetilde W_{\bar \sigma^2_{K_\epsilon t}}|^2\big] 
   = 2t(c_\epsilon t + \bar \sigma_{K_\epsilon}^2).$$
   To handle the probability term $\mathbb{P}_0(t\ge\tau_\epsilon)$,  note that on the event $\{\tau_\epsilon \le t\}$ it holds 
    \begin{align}
       X_{\tau_\epsilon}  &= X_{t \wedge \tau_\epsilon} \le 2 \sup_{s\le t\wedge \tau_\epsilon} \Big| \int_0^s \mu(X_u)\ud u + M_s \Big| \le 2c_\epsilon t+ 2\sup_{u\le \bar \sigma^2_{K_\epsilon}t} |\widetilde W_u|,  
   \end{align}
and by path continuity of $X$ we deduce
   $$\{\tau_\epsilon \le t\} = \{\tau_\epsilon \le t, X_{\tau_\epsilon} = \epsilon\} \ = \Big\{\tau_\epsilon \le t,  \sup_{u \le \bar \sigma^2_{K_\epsilon}t}|\widetilde W_u| \ge \frac{\epsilon}{2} -c_\epsilon t\Big\} \subseteq \Big\{\sup_{u \le \bar \sigma^2_{K_\epsilon}t}|\widetilde W_u| \ge \frac{\epsilon}{2} -c_\epsilon t\Big\}.$$
Without loss of generality, we may choose $t$ sufficiently small so that $\frac{\epsilon}{2} - c_\epsilon t>0$. Hence
$$ \mathbb{P}_0(t\ge\tau_\epsilon) \le \mathbb{P}_0\bigg(\sup_{u \le \bar \sigma^2_{K_\epsilon}t}|\widetilde W_u| \ge\frac{\epsilon}{2} - c_\epsilon t\bigg) \le 4   \mathbb{P}_0 \bigg(\widetilde W_{\bar \sigma^2_{K_\epsilon}t} \ge \frac{\epsilon}{2} - c_\epsilon t\bigg) \le 4 e^{-\frac{(\epsilon -2c_\epsilon t)^2}{8 \bar \sigma^2_{K_\epsilon}t }},$$
where the second to last inequality follows by the symmetry and reflection principle of Brownian motion, while the last one is an application of the so-called Chernoff bound on the Gaussian tail, by which $\mathbb{P}(Z \ge z) \le  e^{-z^2/2}$ for $Z\sim N(0,1)$.  

Finally, we verify the integrability property $\mathbb{E}_0[(L_T^0)^{2p}]^{\frac{1}{p}}<\infty$. By standard inequalities  analogous to those used above
$$ (L_T^0)^{2p} \le 2^{2p-1} \Bigg( \Big(\int_0^T |\mu(X_u)|\ud u \Big)^{2p} +\sup_{s\le T} |M_s|^{2p} \Bigg),$$
hence 
\begin{align*}
    \mathbb{E}_0[(L_T^0)^{2p}] &\le 2^{2p-1} \Bigg( \mathbb{E}_0\Big[\Big(\int_0^T |\mu(X_u)|\ud u \Big)^{2p} \Big] +\mathbb{E}_0\Big[\sup_{s\le T} |M_s|^{2p}\Big] \Bigg)\\
    &\le 2^{2p-1} \Bigg( \mathbb{E}_0\Big[\Big(\int_0^T |\mu(X_u)|\ud u \Big)^{2p} \Big] + C_{2p}\mathbb{E}_0\Big[ \Big(\int_0^T\sigma^2(X_u)\ud u\Big)^p\Big] \Bigg),
\end{align*} 
where the last bound holds by Burkholder–Davis–Gundy inequality applied to the continuous local martingale $M$. Without loss of generality, we may choose $p>1, T>t>0$ according to Assumption $\ref{assumpt: refl_sde}$-$(\ref{assumpt: refl_sde_i})$ and conclude $\mathbb{E}_0[(L_T^0)^{2p}]<\infty$.

Plugging these bounds into \eqref{eq: L_upper}, we get
$$ \mathbb{E}_0[(L_t^0)^2]  \le  2t(c_\epsilon t + \bar \sigma_{K_\epsilon}^2)+ 4^{\frac{1}{q}} \mathbb{E}_0[(L_T^0)^{2p}]^{\frac{1}{p}} e^{-\frac{(\epsilon -2c_\epsilon t)^2}{8 q \bar \sigma^2_{K_\epsilon}t }},$$
leading to 
$$ 0 \le \limsup_{t\downarrow 0} \frac{1}{t} \mathbb{E}_0[(L_t^0)^2] \le 2\bar \sigma^2_{K_\epsilon} + 4^{\frac{1}{q}}\mathbb{E}_0[(L_T^0)^{2p}]^{\frac{1}{p}}  \limsup_{t\downarrow0}   \frac{1}{t} e^{-\frac{(\epsilon -2c_\epsilon t)^2}{8 q \bar \sigma^2_{K_\epsilon}t }}  =  2\bar \sigma^2_{K_\epsilon}.$$
\end{proof}

\section{Proof of the main theorem}\label{sect: proof_main_thm}

This section is devoted to the proof of Theorem~\ref{thm: main}, which proceeds 
in two steps. In Step~(i), we show that the process $\widetilde X$ constructed in 
Section~\ref{sect: setting_main_result} realizes a Feller $(\mu, \sigma)$--reflected elastic diffusion (cf.~Def.~\ref{def: Feller RED}). Since 
$\widetilde X$ has already been shown to be Feller, it remains only to verify that every 
$f \in D(\widetilde{\mathcal{L}})$ satisfies the Robin boundary condition~\eqref{eq: robin_cond} with $\hat S = \widetilde S$, the scale function of $\widetilde X$. 
In Step~(ii), we establish uniqueness (in law). 

\begin{proof}
 \textit{Step (i):} Since $f \in  D(\widetilde{\mathcal{L}})$ is such that $\widetilde{\mathcal{L}} f \in C_0([0,\infty))$ (cf.~\eqref{eq: dom_inf_gen_target}), the differential representation of $\widetilde{\mathcal{L}}$ in \eqref{eq: inf_gen_as_diff_op} yields  $D(\widetilde{\mathcal{L}}) \subseteq 
C^2((0,\infty))$. Furthermore, $\widetilde{\mathcal{L}} f \in C_0([0,\infty))$ implies that $\widetilde{\mathcal{L}} f$ is continuous at the boundary with a finite value $\widetilde{\mathcal{L}} f(0) = \lim_{x \downarrow 0} \widetilde{\mathcal{L}} f(x)$, where 
\begin{equation}\label{eq: st_calc_1}
    \widetilde{\mathcal{L}}f(0) = \lim_{t\downarrow 0} \frac{P_t f(0) - f(0)}{t}
\end{equation}
as implied by the uniform convergence on $[0,\infty)$ of $\frac{P_tf - f}{t}$ to $\widetilde{\mathcal{L}}f$ for $f \in D(\widetilde{\mathcal{L}})$.

We now show that any $f\in D(\widetilde{\mathcal{L}})$ must obey the Robin
condition \eqref{eq: robin_cond} with $\hat S=\widetilde S$. We first prove that $\frac{df}{d\widetilde S}^+(0+)$ exists. Let $f\in
D(\widetilde{\mathcal{L}})$. By \eqref{eq: inf_gen_as_diff_op}, dividing by
$S'$ and using $\sigma>0$ on $[0,\infty)$ (cf.~Assumption \ref{assumpt: refl_sde}),
\[
\frac{d}{dx}\frac{f'(x)}{S'(x)}
= \frac{2\widetilde{\mathcal{L}}f(x)}{\sigma^2(x)S'(x)}
= m(x)\,\widetilde{\mathcal{L}}f(x), \qquad x\in(0,\infty).
\]
Hence, for any fixed $x_0\in(0,\infty)$ and $x\in(0,x_0)$,
\[
f'(x)=S'(x)\left[\frac{f'(x_0)}{S'(x_0)}-\int_x^{x_0}m(z)\,
\widetilde{\mathcal{L}}f(z)\,\mathrm{d}z\right].
\]
Since $\widetilde{\mathcal{L}}f\in C_0(\mathcal{I})\subseteq C_b(\mathcal{I})$ and
$\int_0^{x_0}m(z)\,\mathrm{d}z<\infty$ (as $m\in C([0,\infty))$), the integral
converges as $x\downarrow0$, and  $S'(0+)\in(0,\infty)$ (see Remark~\ref{rem: scale_speed}). Therefore $f'(0+):=\lim_{x\downarrow0}f'(x)$
exists. Furthermore, since  $S=\widetilde S$ on $(0,\infty)$ (see Remark~\ref{rem: tildeX_scale_speed}),
$\widetilde S$ is differentiable with $\widetilde S'=S'$ on $(0,\infty)$;
in particular $\widetilde S'(0+)=S'(0+)\in(0,\infty)$. Hence
\begin{equation}\label{eq: limit_existence}
  \frac{df}{d\widetilde S}^+(0+)=\lim_{x\downarrow 0} \frac{df}{d\widetilde S}^+(x) = \lim_{x\downarrow0}\frac{f'(x)}{\widetilde S'(x)}
=\frac{f'(0+)}{S'(0+)}\in\mathbb{R}.
\end{equation}
Moreover, since $f$ is continuous at $0$, the mean value theorem gives
$\frac{f(h)-f(0)}{h}=f'(\xi_h)$ with $\xi_h\in(0,h)$: since  $\xi_h \downarrow 0$ as $h\downarrow 0$, we have $f'(0)=f'(0+)$,
i.e.\ $f\in C^1([0,\infty))$.

 Assume, for the sake of
contradiction, that
\begin{equation}\label{eq: contr_arg}
\frac{df}{d\widetilde S}^+(0+)=\frac{f'(0)}{S'(0+)}\neq kf(0),
\end{equation}
where $k$ was fixed in the construction of $\widetilde X$ (see
Section~\ref{sect: setting_main_result}). Equivalently,
\begin{equation}\label{eq: contr_arg_rewritten}
f'(0)\neq\lambda f(0),\qquad \lambda=kS'(0+).
\end{equation}

Let 
    \begin{equation}
    \label{eq: loc_seq}
        \tau_n:= \inf\{ t\ge 0: \int_0^t |e^{-\lambda L_s^0} f'(X_s) \sigma(X_s)|^2 \ud s \ge n\} \wedge n.
    \end{equation}
    Since $X$ has continuous paths, on any $[0,t]$, $t\ge 0$ it evolves in a compact set. Recalling $f \in C^1([0,\infty))$ and local boundedness of $\sigma$ (as implied by continuity, cf.~ Assumption \ref{assumpt: refl_sde}) it follows that the integrand in \eqref{eq: loc_seq} is pathwise bounded on any finite interval $[0,t]$. Hence, $\tau_n \uparrow \infty$ a.s. as $n\to \infty$. 

    Setting $X_0= x= 0$, an application of It\^o's formula to $e^{-\lambda L_s^0} f(X_s)$ over $[0, t \wedge \tau_n]$ yields:
\begin{align}\label{eq: ito}
    \mathbb{E}_0 \left[e^{-\lambda  L^0_{t\wedge \tau_n}}f(X_{t\wedge \tau_n})  \right] - f(0) &= \mathbb{E}_0 \left[\int_0^{t\wedge \tau_n}  e^{-\lambda L^0_s}{\mathcal{L}^X}f(X_s)\ud s \right] \\&+ \big(f'(0) -\lambda f(0)\big) \mathbb{E}_0\left[ \int_0^{t\wedge \tau_n} e^{-\lambda  L_s^0} \ud L^0_s \right], 
\end{align}
where $\mathcal{L}^X$ denotes the infinitesimal generator of $X$, such that $ \mathcal{L}^Xf= \widetilde{\mathcal{L}}f$ on $(0,\infty)$ by \eqref{eq: inf_gen_as_diff_op}. Since $f, \widetilde{\mathcal{L}}f \in C_0(\mathcal{I})$ and 
$$\int_0^{t\wedge \tau_n} e^{-\lambda  L_s^0} \ud L^0_s = \frac{1-e^{-\lambda L^0_{t \wedge \tau_n}}}{\lambda} \in [0,1/\lambda),$$ 
in \eqref{eq: ito} we let $n \to \infty$ and apply dominated convergence to obtain 
$$ P_tf(0) - f(0) =  \mathbb{E}_0 \left[\int_0^{t}  e^{-\lambda L^0_s}{\widetilde{\mathcal{L}}}f(X_s)\ud s \right] + \frac{\big(f'(0) -\lambda f(0)\big)}{\lambda} \mathbb{E}_0\left[ 1- e^{-\lambda  L^0_t} \right],$$
hence 
\begin{equation}\label{eq: st_calc_2}
    \frac{P_tf(0) - f(0)}{t} =  \frac{1}{t}\mathbb{E}_0 \left[\int_0^{t}  e^{-\lambda L^0_s}{\widetilde{\mathcal{L}}}f(X_s)\ud s \right] + \frac{\big(f'(0) -\lambda f(0)\big)}{\lambda} \frac{1}{t}\mathbb{E}_0\left[ 1- e^{-\lambda  L^0_t} \right],
\end{equation}
where the term $1 - e^{-\lambda L^0_t}$ appearing in the second expectation can be expanded via a second-order Taylor expansion of $g(y) := e^{-\lambda y}$ on $[0, L^0_t]$ with Lagrange remainder: namely, $g(L^0_t) = g(0) + g'(0) L^0_t + \frac{1}{2}g''(\eta)(L^0_t)^2$ for some $\eta \in (0, L^0_t)$, which yields
\begin{equation}\label{eq: st_calc_3}
    1-e^{-\lambda L^0_t} = \lambda L^0_t - \frac{\lambda^2}{2} e^{-\lambda \eta} (L^0_t)^2 \ge \lambda L^0_t - \frac{\lambda^2}{2} (L^0_t)^2.
\end{equation}

We now separately address the cases  $f'(0) > \lambda f(0)$ and $f'(0) < \lambda f(0)$, exploiting the \emph{local time asymptotics} from Proposition \ref{prop: local_time_asymptotics}.

If $f'(0) >\lambda f(0)$, taking $\liminf$ as $t\downarrow 0$ on both sides of \eqref{eq: st_calc_2}, by \eqref{eq: st_calc_1} and \eqref{eq: st_calc_3} we get
$$ \widetilde{\mathcal{L}}f(0) \ge \widetilde{\mathcal{L}}f(0) + \big(f'(0) -\lambda f(0)\big) \liminf _{t\downarrow 0} \frac{1}{t}\mathbb{E}_0\left[ L^0_t - \frac{\lambda}{2} {(L^0_t)}^2  \right].$$
Using the well-known bound $\liminf_{n}(a_n - b_n) \ge \liminf_{n } a_n -\limsup_{n } b_n $ for real sequences $(a_n)_n$ and $(b_n)_n$, we obtain
$$ 0\ge  \big(f'(0) -\lambda f(0)\big) \bigg( \liminf _{t\downarrow 0} \frac{1}{t}\mathbb{E}_0[L^0_t]  - \frac{\lambda}{2}\limsup_{t\downarrow 0} \frac{1}{t} \mathbb{E}_0\big[(L^0_t)^2 \big] \bigg), $$
which, by \eqref{eq: liminf_bound} and \eqref{eq: limsup_bound}, yields the contradiction $0\ge \infty$. 
 
Instead, if $f'(0) <\lambda f(0)$, by the same arguments we obtain
$$ 0 \le  \big(f'(0) -\lambda f(0)\big) \liminf _{t\downarrow 0} \frac{1}{t}\mathbb{E}_0\left[ L^0_t - \frac{\lambda}{2} {(L^0_t)}^2  \right] \le  \big(f'(0) -\lambda f(0)\big)  \bigg( \liminf _{t\downarrow 0} \frac{1}{t}\mathbb{E}_0[L^0_t]  - \frac{\lambda}{2}\limsup_{t\downarrow 0} \frac{1}{t} \mathbb{E}_0\big[(L^0_t)^2 \big] \bigg), $$
leading, by \eqref{eq: liminf_bound} and \eqref{eq: limsup_bound}, to the contradiction $0\le -\infty$.
We conclude that any $f\in D(\widetilde{\mathcal{L}})$ satisfies the Robin condition as in \eqref{eq: robin_cond} with $\hat S = \widetilde S$. Hence $\widetilde X$ realizes a Feller $(\mu, \sigma)$--RED. 
\medskip 

\noindent \textit{Step (ii):}
By virtue of \textit{Step (i)}, it holds that
\begin{align*}
    D(\widetilde{\mathcal{L}}) &= \Big\{ f \in C_0(\mathcal{I}): \exists g \in C_0(\mathcal{I}) \text{ s.t. } \lim_{t\downarrow 0}\Big|\Big| \frac{P_t f - f}{t} - g\Big|\Big|_\infty = 0 \Big\} \\
    &= \Big\{ f \in C_0(\mathcal{I}) \cap C^2((0,\infty)): \exists g \in C_0(\mathcal{I}) \text{ s.t. } \lim_{t\downarrow 0}\Big|\Big| \frac{P_t f - f}{t} - g\Big|\Big|_\infty = 0, ~ \frac{df}{d\widetilde S}^+(0+) = kf(0) \Big\}
\end{align*}  
where we recall that the first equality holds by definition of the infinitesimal generator of a Feller process (see \eqref{eq: dom_inf_gen_target}).  

Suppose there exists another Feller $(\mu, \sigma)$--RED (in the sense of Definition ~\ref{def: Feller RED}) on $\mathcal{I}$, denoted by $X^*$,  sharing the same interior dynamics as $\widetilde X$, that is, on $(0,\infty)$ the process $X^*$ is governed by the same infinitesimal coefficients $(\mu, \sigma^2)$ characterizing $\widetilde X$ (cf.~ Def.~\ref{def: weak_sol_refl_sde}) and satisfies Assumption \ref{assumpt: refl_sde}. As such, on $(0,\infty)$ the scale and speed density of $X^*$ coincide with $S'$ and $m$ (respectively) as defined in \eqref{eq: scale_speed_X}, and satisfy the properties of Remark \ref{rem: scale_speed} (see also Remark \ref{rem: tildeX_scale_speed}). Let $(P^*_t)_{t\ge 0}$ be the (Feller) semigroup of $X^*$ and $(\mathcal{L^*}, D(\mathcal{L}^*))$ its infinitesimal generator (cf.~\eqref{eq: dom_inf_gen_target}-\eqref{eq: inf_gen_target} wih $E= \mathcal{I}$). By the same arguments adopted for the process $\widetilde X$ at the end of Section \ref{sect: tildeX_properties}, on $(0,\infty)$ the operator $\mathcal{L}^*$ coincides with the second-order differential operator
$$ \mathcal{L}^*f (x) = \frac{\sigma^2(x)}{2} f''(x) + \mu(x) f'(x), \quad x \in (0,\infty), f \in \mathcal{D(\mathcal{L}^*})$$
implying $D(\mathcal{L}^*) \subset C^2((0,\infty))$. Moreover, denoting by $S^*$ the scale function of $X^*$, any $f \in D(\mathcal{L}^*)$ satisfies $\frac{df}{dS^*}^+(0+)= kf(0)$ by definition of RED. It follows that 
\begin{align*}
    D(\mathcal{L}^*) = \Big\{ f \in C_0(\mathcal{I}) \cap C^2((0,\infty)): \exists g \in C_0(\mathcal{I}) \text{ s.t. } \lim_{t\downarrow 0}\Big|\Big| \frac{P^*_t f - f}{t} - g\Big|\Big|_\infty = 0, ~ \frac{df}{dS^*}^+(0+) = kf(0) \Big\}.
\end{align*}  
We now aim to show $(\widetilde{\mathcal{L}},D(\widetilde{\mathcal{L}})) = (\mathcal{L}^*, D(\mathcal{L}^*))$, that is, $D(\widetilde{\mathcal{L}}) = D(\mathcal{L}^*)$ and $\widetilde{\mathcal{L}} f = \mathcal{L}^*f$ for all $f \in D(\widetilde{\mathcal{L}}) = D(\mathcal{L}^*)$. 
For any $g \in C_0(\mathcal{I})$ and $\alpha > 0$, consider the boundary value problem (BVP)
\begin{equation}\label{eq: bvp}
\begin{cases} 
\displaystyle\alpha f(x) - \frac{\sigma^2(x)}{2}f''(x) - \mu(x)f'(x) = g(x), & x \in (0,\infty), \\[10pt] 
k f(0) = \displaystyle\lim_{x\downarrow 0}\frac{f'(x)}{S'(x)}, \\[10pt] 
f \in C_0(\mathcal{I}) \cap C^2((0,\infty)),
\end{cases}
\end{equation}
where we recall that $S'$ denotes the scale density of $X$ (cf.~\eqref{eq: scale_speed_X}). The first equation can be equivalently written as $(\alpha - \widetilde{\mathcal{L}}) f = g$ and $(\alpha - \mathcal{L}^*) f = g$.

Since $\widetilde{\mathcal{L}}$ and $\mathcal{L}^*$ are the infinitesimal generators of contraction semigroups (cf.~\cite[Ch.~3, Def.~2.1]{revuz2013continuous}), the Hille-Yosida theorem (see, e.g., \cite[Ch.~1, Thm.~3.1]{Pazy1992SemigroupsOL}) guarantees that $(0,\infty)$ is contained in their resolvent sets, that is, for any fixed $\alpha>0$, $\widetilde{R}_\alpha := (\alpha - \widetilde{\mathcal{L}})^{-1}$ and $R^*_\alpha := (\alpha - \mathcal{L}^*)^{-1}$ are well-defined operators on $C_0(\mathcal{I})$, mapping into $D(\widetilde{\mathcal{L}})$ and $D(\mathcal{L}^*)$, respectively. Hence $\widetilde{R}_\alpha g$ and $R^*_\alpha g$ both solve the first equation on $(0,\infty)$ and satisfy the regularity requirement. They also satisfy the Robin condition: every $f \in D(\widetilde{\mathcal{L}})$ (resp.\ $D(\mathcal{L}^*)$) satisfies \eqref{eq: robin_cond} with $\hat S = \widetilde S$ ($\hat S = S^*$), and
\[
\frac{df}{d\widetilde{S}}^+(0+) = \lim_{x \downarrow 0} \frac{f'(x)}{\widetilde{S}'(x)} = \lim_{x \downarrow 0} \frac{f'(x)}{S'(x)},
\]
the limit existing by the arguments preceding \eqref{eq: limit_existence}, which apply identically to $\frac{df}{d{S^*}}^+$. Thus, both $\widetilde R_\alpha g$ and $R^*_\alpha g$ solve the BVP \eqref{eq: bvp} for any $g \in C_0(\mathcal{I})$ and $\alpha > 0$. We wish to conclude $\widetilde R_\alpha g = R^*_\alpha g$. 

To this end, the uniqueness of the solution to \eqref{eq: bvp} follows from standard arguments: let $h_1, h_2$ be two solutions to \eqref{eq: bvp} and let $h := h_1 - h_2$. It suffices to show $h \equiv 0$. By linearity, $h$ also solves \eqref{eq: bvp} with $g \equiv 0$. Suppose $\sup_{x \in \mathcal{I}} h(x) > 0$: since $h \in C_0(\mathcal{I})$, this supremum is attained at some $x_0 \in \mathcal{I}$. 
If $x_0 \in (0,\infty)$, then $h'(x_0) = 0$ and $h''(x_0) \le 0$, so that $\alpha h(x_0) = \frac{\sigma^2(x_0)}{2}h''(x_0) + \mu(x_0)h'(x_0) \le 0$, contradicting $h(x_0) > 0$ and $\alpha > 0$. 
If $x_0 = 0$, then $h(x) \le h(0)$ for all $x>0$, hence
$\frac{h(x) - h(0)}{S(x) - S(0)}\le 0$. 
By the (extended) mean value theorem, for each $x>0$ there exists $\xi_x\in (0,x)$ such that
$\frac{h(x) - h(0)}{S(x) - S(0)} = \frac{h'(\xi_x)}{S'(\xi_x)}$.
Since $\xi_x \downarrow 0$ as $x\downarrow 0$ and $\lim_{y\downarrow0}\frac{h'(y)}{S'(y)}=kh(0)$, we have
$\lim_{x\downarrow0}\frac{h'(\xi_x)}{S'(\xi_x)}=kh(0)>0$, contradicting
$\frac{h(x) - h(0)}{S(x) - S(0)}\le 0$. Hence $h \le 0$, and by symmetry $h \ge 0$, so $h \equiv 0$.
We conclude $\widetilde R_\alpha  = R^*_\alpha$ on $C_0(\mathcal{I})$.

Since the maps $ (\alpha - \widetilde{\mathcal{L}}): D(\widetilde{\mathcal{L}})\to C_0(\mathcal{I})$ and $(\alpha - {\mathcal{L}}^*): D({\mathcal{L}}^*)\to C_0(\mathcal{I}) $ are \emph{bijections} (cf.~\cite[Ch.~7, Prop.~1.4]{revuz2013continuous}), we have $\text{Range}(\widetilde R_\alpha) = D(\widetilde{\mathcal
L})$ and $\text{Range}( R^*_\alpha) = D(\mathcal
L^*)$, thus  $\widetilde R_\alpha  = R^*_\alpha$ implies $D(\widetilde{\mathcal
L}) =D(\mathcal
L^*)$. Finally, the equality between $\widetilde{\mathcal
L}$ and $\mathcal{L}^*$ on the open interval $(0,\infty)$ is evident, and since they both map into $C_0(\mathcal{I})$, the equality extends at the boundary by continuity. It follows  $(\widetilde{\mathcal{L}},D(\widetilde{\mathcal{L}}) ) = (\mathcal{L}^*, D(\mathcal{L}^*))$.

It is known that a strongly continuous semigroup of bounded linear operators is uniquely 
determined by its infinitesimal generator (see, e.g., \cite[Ch.1, Thm. 2.6]{Pazy1992SemigroupsOL}). 
Since $\widetilde{X}$ and $X^*$ are Feller processes, their semigroups 
$(\widetilde{P}_t)_{t\geq 0}$ and $(P^*_t)_{t\geq 0}$ are both strongly continuous 
contractions on $C_0(\mathcal{I})$, and we have shown their infinitesimal 
generator to coincide.  We deduce $\widetilde{P}_t = P^*_t$ for all $t \geq 0$. 
Moreover, as recalled at the beginning of Section \ref{sect: setting_main_result},  each Feller semigroup uniquely identifies a homogeneous transition function 
(cf.\ \cite[Ch.~3, Prop.~2.2]{revuz2013continuous}), implying that $\widetilde{X}$ and $X^*$ share 
the same transition function. As a Markov process is uniquely determined in law by 
its initial distribution and its transition function 
(cf.\ \cite[Ch.~3, Thm.~1.5]{revuz2013continuous}),  for each fixed starting 
point $x \in \mathcal{I}$ we conclude $\widetilde{X} \stackrel{d}{=} X^*$.
\end{proof}

\medskip
\paragraph{\textbf{Acknowledgements.}} The author is grateful to T.~De Angelis
and A.~Milazzo for their guidance and constant feedback during the
development of this work. The author also acknowledges financial
support from the PhD fellowship of the University of Turin.

\bibliographystyle{abbrv}
\bibliography{Bibliography}
\end{document}